\documentclass{amsart}
\usepackage{amsmath,amsthm,amssymb,amsfonts,amscd}
\input xy
\xyoption{all}

\theoremstyle{plain}
\newtheorem{thm}{Theorem}
\newtheorem{prp}[thm]{Proposition}

\theoremstyle{definition}
\newtheorem{dfn}[thm]{Definition}

\theoremstyle{remark}
\newtheorem*{rmk}{Remark}

\newtheorem*{ack}{Acknowledgements}

\newcommand{\Hom}{\textup{Hom}}
\newcommand{\End}{\textup{End}}
\newcommand{\Aut}{\textup{Aut}}
\newcommand{\Inn}{\textup{Inn}}

\newcommand{\id}{\textup{id}}
\newcommand{\Img}{\textup{Im}}
\newcommand{\Ker}{\textup{Ker}}
\newcommand{\GL}{\textup{GL}}
\newcommand{\PGL}{\textup{PGL}}

\newcommand{\ch}{\textup{char }}

\newcommand{\SD}{\textup{SD}}

\newcommand{\F}{\mathcal{F}}

\begin{document}
\title[The weighted fusion category algebra and the $q$-Schur algebra]{The weighted fusion category algebra and the $q$-Schur algebra for $\GL_2(q)$}
\author{Sejong Park}
\address{University of Aberdeen, Aberdeen, UK AB24 3UE}
\email{spark@maths.abdn.ac.uk}

\begin{abstract}
We show that the weighted fusion category algebra of the principal $2$-block $b_0$ of $\GL_2(q)$ is the quotient of the $q$-Schur algebra $\mathcal{S}_2(q)$ by its socle, for $q$ an odd prime power.  As a consequence, we get a canonical bijection between the set of isomorphism classes of simple $k\GL_2(q)b_0$-modules and the set of conjugacy classes of $b_0$-weights in this case.
\end{abstract}
\maketitle

\section{Introduction} \label{S:Intro}
Let $k$ be an algebraically closed field of characteristic $l>0$. Let $q$ be a prime power which is coprime to $l$.  Consider two $k$-algebras associated with $\GL_n(q)$, namely the weighted fusion category algebra $\overline{\mathcal{F}}(b_0)$ of the principal $l$-block $b_0$ of $\GL_n(q)$ defined by Linckelmann~\cite{Linckelmann2004} and the $q$-Schur algebra $\mathcal{S}_n(q)$ introduced by Dipper and James~\cite{DipperJames1989}. Since both are quasi-hereditary and carry informations of representations of $\GL_n(q)$(Theorems \ref{T:FusionProperty} and \ref{T:SchurProperty}), one may conjecture that there is a certain relation between them.

We give definitions and some properties of these algebras in Sections~\ref{S:Fusion} and \ref{S:Schur}.  In Section~\ref{S:Example}, we compute the Morita types of these algebras in a special case and find a relation between them: 

\begin{thm} \label{T:Main}
Let $k$ be an algebraically closed field of characteristic $2$ and let $q$ be an odd prime power. Then the weighted fusion category algebra $\overline{\mathcal{F}}(b_0)$ over $k$ of the principal $2$-block $b_0$ of $\GL_2(q)$ is Morita equivalent to the quotient of the $q$-Schur algebra $\mathcal{S}_2(q)$ over $k$ by its socle.
\end{thm}

Theorem~\ref{T:Main} implies a canonical bijection between the set of isomorphism classes of simple $k\GL_2(q)b_0$-modules and the set of isomorphism classes of simple $\overline{\F}(b_0)$-modules, which in turn is in a bijective correspondence with the set of conjugacy classes of $b_0$-weights.  We discuss this canonical bijection in more detail in Section~\ref{S:RemarkNaturalBijection}.

\section{The weighted fusion category algebra} \label{S:Fusion}
We summarize the construction of the weighted fusion category algebra and its basic properties,  following Linckelmann \cite{Linckelmann2004}.  We restrict our attention to the principal block case and avoid discussing the ``twisting and gluing'' procedure.(See \cite[{4.1--4.4, 5.1}]{Linckelmann2004})

Let $k$ be an algebraically closed field of characteristic $l>0$ and let $G$ be a finite group. Let $b_0$ be the principal $l$-block of $G$, i.e.\ the unique primitive idempotent in the center $Z(kG)$ of the group algebra $kG$ which is not contained in the augmentation ideal of $kG$. Fix a defect group $P$ of $b_0$, namely a Sylow $l$-subgroup of $G$. The \textit{fusion system} of the block $b_0$ (on $P$) is the same as the fusion system of the group $G$ (on $P$): it is the category $\F = \F_P(G)$ whose objects are the subgroups of $P$ and such that for each pair $Q, R$ of subgroups of $P$ the morphism set $\Hom_{\mathcal{F}}(Q,R)$ consists of the group homomorphisms from $Q$ to $R$ induced by conjugations in $G$. It is independent of the choice of a defect group $P$, up to equivalence of categories.

The \textit{orbit category} of $\mathcal{F}$ is the category $\overline{\mathcal{F}}$ whose objects are again the subgroups of $P$ and such that for each pair $Q, R$ of subgroups of $P$ the morphism set $\Hom_{\overline{\mathcal{F}}}(Q,R)$ consists of the orbits of the group of inner automorphisms $\Inn(R)$ of $R$ in $\Hom_{\mathcal{F}}(Q,R)$. 

A subgroup $Q$ of $P$ is called $\mathcal{F}$\textit{-centric} if every subgroup $R$ of $P$ which is  $\mathcal{F}$-isomorphic to $Q$ is centric in $P$, i.e.\ $C_P(R) = Z(R)$. We denote by $\overline{\mathcal{F}}^c$ the full subcategory of the orbit category $\overline{\mathcal{F}}$ consisting of $\mathcal{F}$-centric subgroups of $P$.

Let $k\overline{\mathcal{F}}^c$ be the \textit{category algebra} of $\overline{\mathcal{F}}^c$ over $k$, that is, the $k$-algebra whose $k$-basis consists of morphisms of $\overline{\mathcal{F}}^c$ and such that multiplication is induced by composition of morphisms.

\begin{dfn} \label{D:Fusion}
With the above notations, let $\overline{e} = \sum_{Q} \overline{e}_Q$ where $Q$ runs over all $\mathcal{F}$-centric subgroups of $P$ and $\overline{e}_Q$ denotes the sum of all defect zero blocks of $k \Aut_{\overline{\mathcal{F}}}(Q)$. Then the \textit{weighted fusion category algebra} of the block $b_0$ is the truncated algebra
\[
	\overline{\mathcal{F}}(b_0) = \overline{e} k\overline{\mathcal{F}}^c \overline{e}.
\]
\end{dfn}

The significance of the weighted fusion category algebra is summarized in the following theorem:

\begin{thm}[{\cite[4.5, 5.1]{Linckelmann2004}}] \label{T:FusionProperty}
Let $k$ be an algebraically closed field of characteristic $l > 0$ and let $b_0$ be the principal $l$-block of a finite group $G$.  Then the weighted fusion category algebra $\overline{\mathcal{F}}(b_0)$ over $k$ of the block $b_0$ is quasi-hereditary and Alperin's weight conjecture for the block $b_0$ is equivalent to the equality
\[
	l(kGb_0) = l(\overline{\mathcal{F}}(b_0))
\]
where $l(A)$ denotes the number of isomorphism classes of simple $A$-modules for a finite dimensional $k$-algebra $A$.
\end{thm}

\begin{rmk}
We refer to Alperin's original paper~\cite{Alperin1987} for the definition of weights and the statement of Alperin's weight conjecture.  Note that, for an $\F$-centric subgroup $Q$ of $P$, the defect zero blocks of $k\Aut_{\overline{\F}}(Q)$ appearing in Definition~\ref{D:Fusion} (if any) correspond to the $b_0$-weights having $Q$ as their first component.  Alperin's weight conjecture is positively confirmed, among others, for finite general linear groups by Alperin and Fong~\cite{AlperinFong1990} in odd characteristics and by An~\cite{An1992} in characteristic $2$.
\end{rmk}

\section{The $q$-Schur algebra} \label{S:Schur}
We review the definition and some basic properties of the $q$-Schur algebra defined by Dipper and James \cite{DipperJames1989}, following the presentation of Mathas \cite{Mathas1999}.

Let $k$ be a field and let $q$ be a nonzero element of $k$. The \textit{Iwahori-Hecke algebra} of the symmetric group $\Sigma_n$ on $n$ letters is the $k$-algebra $\mathcal{H} = \mathcal{H}_{k,q}(\Sigma_n)$ whose $k$-basis is $\{ T_w \mid w \in \Sigma_n \}$ and such that multiplication is given by
\begin{equation*}
	T_w T_s =
	\begin{cases}
		T_{ws},		&\text{if $l(ws) > l(w)$,}\\
		q T_{ws} + (q-1) T_w,	&\text{if $l(ws) < l(w)$,}
	\end{cases}
\end{equation*}
where $w \in \Sigma_n$, $s = (i, i+1) \in \Sigma_n$ for some $0 < i < n$, and $l(w)$ is the length of $w$. 

A \textit{composition} of $n$ is a sequence $\mu = (\mu_1, \mu_2, \dots)$ of nonnegative integers $\mu_i$ whose sum is equal to $n$. The $height$ of a composition $\mu$ is the smallest positive integer $d$ such that $\mu_{d+1} = \mu_{d+2} = \dots = 0$. For a composition $\mu$ of $n$ with height $d$, let $\Sigma_{\mu}$ be the corresponding \textit{Young subgroup} of $\Sigma_n$ isomorphic to $\Sigma_{\mu_1} \times \Sigma_{\mu_2} \times \dots \times \Sigma_{\mu_d}$. Set $m_{\mu} = \sum_{w \in \Sigma_{\mu}} T_w$ and let $M^{\mu} = m_{\mu} \mathcal{H}$, the right $\mathcal{H}$-submodule of $\mathcal{H}$ generated by $m_{\mu}$.

\begin{dfn} \label{D:Schur}
Let $\Lambda(d,n)$ be the set of all compositions of $n$ with height $\leq d$. Then the $q$\textit{-Schur algebra} is the endomorphism algebra
\[
	\mathcal{S}_{d,n}(q) = \End_{\mathcal{H}} \Bigl( \bigoplus_{\mu \in \Lambda(d,n)} M^{\mu} \Bigr).
\]
We write $\mathcal{S}_n(q) = \mathcal{S}_{n,n}(q)$.
\end{dfn}

The $q$-Schur algebra has the following properties:

\begin{thm}[{\cite[4.16, 6.47]{Mathas1999}}] \label{T:SchurProperty}
Let $k$ be a field and let $q$ be a nonzero element of $k$. Then the $q$-Schur algebra $\mathcal{S}_{d,n}(q)$ over $k$ is quasi-hereditary. If $\ch k = l > 0$ and $q$ is a prime power which is coprime to $l$, then the decomposition matrix of $k\GL_n(q)$ is completely determined by the decomposition matrices of the $q^r$-Schur algebras $\mathcal{S}_{m}(q^r)$ over $k$ for $rm \leq n$.
\end{thm}

Gruber and Hiss \cite{GruberHiss1997} and Takeuchi \cite{Takeuchi1996} give an alternative way of computing the Morita types of the $q$-Schur algebras. 

\begin{thm}[\cite{GruberHiss1997}, \cite{Takeuchi1996}] \label{T:SchurAlternative}
Let $k$ be an algebraically closed field of characteristic $l > 0$ and let $q$ be a prime power which is coprime to $l$. Let $G = \GL_n(q)$ and let $B$ be the set of all upper triangular matrices in $G$. Then the $q$-Schur algebra $\mathcal{S}_{n}(q)$ over $k$ is Morita equivalent to the image of the $k$-algebra homomorphism
\[
	kG \rightarrow \End_{k} (k[G/B])
\]
sending $a \in kG$ to the $k$-linear endomorphism of $k[G/B]$ given by left multiplication by $a$ on $k[G/B]$.
\end{thm}

\section{The case $G = \GL_2(q)$, $q$ odd, in characteristic $2$} \label{S:Example}

Let $k$ be an algebraically closed field of characteristic $2$ and let $q$ be an odd prime power. In this case, we have

\begin{prp} \label{T:FusionExample}
The weighted fusion category algebra $\overline{\mathcal{F}}(b_0)$ over $k$ of the principal $2$-block $b_0$ of $\GL_2(q)$ is Morita equivalent to the path algebra of the quiver 
\[
\xymatrix@C=20mm{
	^{1}\bullet & \bullet^{2} \ar@<.35ex>[l]_{\alpha} }
\]
\end{prp}

\begin{prp}[{\cite[3.3(A)]{ErdmannNakano2001}}] \label{T:SchurExample}
The $q$-Schur algebra $\mathcal{S}_2(q)$ over $k$ is Morita equivalent to the path algebra of the quiver 
\[
\xymatrix@C=20mm{
	^{1}\bullet \ar@<.05ex>[r]^{\beta} & 
	\bullet^{2} \ar@<.75ex>[l]^{\gamma} }
\]
with relation $\beta\gamma = 0$.
\end{prp}

A proof of Proposition~\ref{T:FusionExample} is given in Sections \ref{S:FusionProof3mod4} and \ref{S:FusionProof1mod4}.  Proposition~\ref{T:SchurExample} is a consequence of more general results of Erdmann and Nakano~\cite{ErdmannNakano2001}.  For the convenience of the reader, we sketch a proof of Proposition~\ref{T:SchurExample} in Section \ref{S:SchurProof}.  Theorem \ref{T:Main} follows immediately from Propositions~\ref{T:FusionExample} and~\ref{T:SchurExample}.

\subsection{Proof of Proposition \ref{T:FusionExample} when $q \equiv 3 \mod 4$} \label{S:FusionProof3mod4}
Let $G=\GL_2(q)$ where $q$ is a prime power such that $q \equiv 3 \mod 4$.  Let $2^{m - 2}$ be the highest $2$-power dividing $q + 1$ and let $\xi$ be a primitive $2^{m - 1}$th root of unity in $\mathbb{F}_{q^2}$. Note that $m \geq 4$. Then the subgroup $P$ of $G$ generated by
\[
	x = 
	\begin{pmatrix}
		0 & 1 \\
		1 & a 
	\end{pmatrix},
	\qquad
	t =
	\begin{pmatrix}
		1 & a \\
		0 & -1
	\end{pmatrix}
	\qquad
	(a = \xi + \xi ^q)
\]
is a Sylow $2$-subgroup of $G$. One immediately checks that $x$ and $t$ are of order $2^{m-1}$ and $2$, respectively, and
\[
	txt = x^{2^{m-2} -1}.
\]
In other words, $P$ is the semidiheral group $\SD_{2^m}$ of order $2^m$.

Let $\F = \F_P(G)$.  Then the $\mathcal{F}$-centric subgroups of $P$ are as follows:
\begin{enumerate}
	\item $C_2 \times C_2 \cong \langle x^{2^{m-2}}, tx^{2i} \rangle$
	\item $D_{2^k} \cong \langle x^{2^{m-k}}, tx^{2i} \rangle \text{ where } 3 \leq k \leq m-1$
	\item $Q_{2^k} \cong \langle x^{2^{m-k}}, tx^{2i+1} \rangle \text{ where } 3 \leq k \leq m-1$
	\item $C_{2^{m-1}} \cong \langle x \rangle$
	\item $P$
\end{enumerate}
Recall that the automorphism groups of cyclic, dihedral, semidihedral, and (generalized) quaternion 2-groups of order $\geq 4$ are all nontrivial 2-groups except for 
\[
	\Aut(C_2 \times C_2) \cong \Sigma_3,\qquad \Aut(Q_8) \cong \Sigma_4.
\]
So the $\overline{\F}$-automorphism group of an $\F$-centric subgroup $R$ of $P$ of type (2), (3) with $k>3$, (4), or (5) is a (possibly trivial) 2-group.  If $R = P$, then since $\Aut_{\overline{\F}}(P)$ is also a $2'$-group, we have $\Aut_{\overline{\F}}(P) = \{\id_{P}\}$ and hence $\overline{e}_P = 1$.  Let $e_1 = \overline{e}_P$.  If $R < P$, then we have $\Inn(R) < \Aut_{\F}(R)$ and hence $\Aut_{\overline{\F}}(R)$ is a nontrivial 2-group.  Therefore $\overline{e}_R = 0$.

Also, since $x^{2^{m-2}} = \left( \begin{smallmatrix} -1 & 0\\ 0 & -1 \end{smallmatrix} \right) \in Z(G)$ and $tx^{2i}$, $-tx^{2i}$ are $G$-conjugate, the $\mathcal{\overline{F}}$-automorphism group of a Klein four subgroup $R$ of $P$ is isomorphic to $C_2$, yielding $\overline{e}_R = 0$. Thus it remains to consider the quaternion subgroups of order $8$. Set
\[
	Q_{i} = \langle x^{2^{m-3}}, tx^{2i+1} \rangle,\quad i = 0,1,\ldots,2^{m-4}-1.
\]
First observe that all $Q_i$ are $P$-conjugate. Indeed, for each pair of indices $i, j$, let $k = ( 2^{m-3} - 1) ( j - i )$. Then
\[
	x^k tx^{2i+1} x^{-k} = t x^{(2^{m-2}-1) k} x^{2i+1-k} = tx^{2j+1}.
\]
So it suffices to consider only $Q:= Q_0 = \langle x^{2^{m-3}}, tx \rangle$. We have $\Aut_{P}(Q) \cong D_8$ and $\Aut(Q) \cong \Sigma_4$. Thus $\Aut_{\mathcal{F}}(Q)$ is either $\Aut_{P}(Q)$ or $\Aut(Q)$. Since $x^{2^{m-3}}$ and $tx$ are $G$-conjugate but automorphisms in $\Aut_{P}(Q)$ do not send $x^{2^{m-3}}$ to $tx$, one finds that 
\[
	\Aut_{\mathcal{F}}(Q) = \Aut(Q) \cong \Sigma_4.
\]
(Note that this also follows from the Frobenius normal $p$-complement theorem.)  Now $\Aut_{\mathcal{\overline{F}}}(Q) = \Aut_{\mathcal{F}}(Q) / \Aut_{Q}(Q)$ and $\Aut_{Q}(Q) \cong C_2 \times C_2$. Thus we have
\[
	\Aut_{\mathcal{\overline{F}}}(Q) \cong \Sigma_3.
\]
Since $k\Sigma_3 \cong kC_2 \times M_2(k)$ as $k$-algebras, one finds that $\overline{e}_{Q} k \Aut_{\mathcal{\overline{F}}}(Q) \overline{e}_{Q} \cong M_2(k)$. Let $e_2$ be the element of $\overline{e}_{Q} k \Aut_{\mathcal{\overline{F}}}(Q) \overline{e}_{Q}$ which corresponds to
$\left(
\begin{smallmatrix}
1 & 0\\
0 & 0
\end{smallmatrix}
\right)$
via this isomorphism.

Set $A := \overline{e} k \mathcal{\overline{F}}^{c} \overline{e}$ where $\overline{e} = \overline{e}_P + \overline{e}_Q$.  Then $A$ is Morita equivalent to $\overline{\mathcal{F}}(b_0)$.  We have a decomposition as $k$-vector spaces
\[
	A = A_1 \oplus A_2 \oplus J 
\]
where
\begin{gather*}
	A_1 = \overline{e}_P k \Aut_{\mathcal{\overline{F}}}(P) \overline{e}_P \cong k,\\
	A_2 = \overline{e}_{Q} k \Aut_{\mathcal{\overline{F}}}(Q) \overline{e}_{Q}
	\cong M_{2}(k),\\
	J = \overline{e}_{P} k \Hom_{\mathcal{\overline{F}}}(Q,P) \overline{e}_{Q}.
\end{gather*}
Since $J^2 = 0$ and $A/J \cong k \times M_{2}(k)$, $J$ is the Jacobson radical of $A$ and there are exactly two nonisomorphic simple $A$-modules $S_1$, $S_2$ with corresponding projective indecomposable $A$-modules $P_1 = A e_1, P_2 = A e_2$. Note that $P_1 = S_1 \cong k$ and $J P_2 / J^2 P_2 = k \Hom_{\mathcal{\overline{F}}} (Q, P) e_2 \cong S_1$. Therefore we get the desired result.

\subsection{Proof of Proposition \ref{T:FusionExample} when $q \equiv 1 \mod 4$} \label{S:FusionProof1mod4}

Let $G = \GL_2(q)$ where $q$ is a prime power such that $q \equiv 1 \mod 4$.  Let $2^{m}$ be the highest $2$-power dividing $q - 1$ and let $\eta$ be a primitive $2^{m}$th root of unity in $\mathbb{F}_{q}$. Note that $m \geq 2$. Then the subgroup $P$ of $G$ generated by
\[
	x = 
	\begin{pmatrix}
		\eta & 0 \\
		0 & 1 
	\end{pmatrix},
	\quad
	y=
	\begin{pmatrix}
		1 & 0 \\
		0 & \eta
	\end{pmatrix},
	\quad
	t =
	\begin{pmatrix}
		0 & 1 \\
		1 & 0
	\end{pmatrix}
\]
is a Sylow $2$-subgroup of $G$. Since $x$, $y$ commute and $txt = y$, we see that $P \cong C_{2^{(m-1)/2}} \wr \Sigma_2$.  Note that $Z_0 := Z(P) = Z(G) \cap P = \langle xy \rangle \cong C_{2^m}$.

Let $\F = \F_P(G)$.  Then the $\mathcal{F}$-centric subgroups of $P$ are as follows:
\begin{enumerate}
	\item $\langle x,y \rangle$
	\item $\langle xy, tx^{i} \rangle \text{ where } \eta^{i} \neq \eta^{2j} \text{ for any integer } j$
	\item $\langle xy, x^{2^i}, tx^{j} \rangle \text{ where } 0 \leq i \leq m-1, 0 \leq j < 2^{i}$
\end{enumerate}
	
Let $R$ be an $\mathcal{F}$-centric subgroup of $P$. If $R = \langle x, y \rangle$, then we have
\[
	\Aut_{\overline{\F}}(R) \cong N_G(R) / RC_G(R) = L\Sigma_2 / L \cong \Sigma_2,
\]
where $L$ denotes the diagonal subgroup of $G$ and $\Sigma_2$ is viewed as the subgroup of the permutation matrices in $G$.

Now suppose that $R$ is of type (2) or (3).  Since $Z_0 \subseteq Z(G)$, elements of $Z_0$ are fixed by any $\mathcal{F}$-morphism. So every $\mathcal{F}$-automorphism of $R$ induces an automorphism of $R/Z_0$, giving rise to a surjective group homomorphism
\[
	\Phi : \Aut_{\mathcal{F}}(R) \twoheadrightarrow \Aut_{G/Z_0}(R/Z_0).
\]
Note that the kernel $\Ker(\Phi)$ of $\Phi$ is isomorphic to a certain subgroup of the group $\Hom(R, Z_0)$ whose multiplication is given by pointwise multiplication.  In particular $\Ker(\Phi)$ is an abelian $2$-group.

If $R$ is of type (2), then $R/Z_0 \cong C_2$, so $\Aut(R/Z_0) = \{\id_{R/Z_0}\}$.  One can easily check that $\Ker(\Phi) \cong C_2$ in this case.  Since $R$ is abelian, it follows that $\Aut_{\overline{\mathcal{F}}}(R) \cong C_2$.

Suppose that $R$ is of type (3). Then $R/Z_0$ is a dihedral $2$-group of order $\geq 4$; it is of order $4$ (i.e.\ a Klein four group) if and only if $i = m-1$. So if $i \neq m-1$, then $R/Z_0$ is a dihedral $2$-group of order $\geq 8$, and hence its automorphism group is a (nontrivial) $2$-group.  Thus $\Aut_{\mathcal{F}}(R)$ is a $2$-group. Now if $R < P$, then $\Inn(R) < \Aut_{\F}(R)$, so $\Aut_{\overline{\mathcal{F}}}(R)$ is a nontrivial $2$-group; if $R = P$, then $\Aut_{\overline{\F}}(P)$ is also a $2'$-group and hence $\Aut_{\overline{\mathcal{F}}}(P) = \{\id_P\}$.

Finally, let $R$ be of type (3) with $i = m-1$. There are two $P$-conjugacy classes among these $\mathcal{F}$-centric subgroups. Indeed, for any $j$, 
\[
	\langle xy, x^{2^{m-1}}, tx^{j} \rangle \cong \langle xy, x^{2^{m-1}}, tx^{j+2} \rangle
\]
because $x^{-1} ( tx^{j+1}y ) x = tx^{j+2}$. Set
\[
	R_1 = \langle xy, x^{2^{m-1}}, t \rangle, \qquad
	R_2 = \langle xy, x^{2^{m-1}}, tx \rangle.
\]
Since $R_{i}/Z_0 (i=1,2)$ is a Klein four group, its full automorphism group is isomorphic to $\Sigma_3$, permuting its three nonidentity elements. Those three nonidentity elements of $R_1/Z_0$ are all $G$-conjugate; in $R_2/Z_0$, the elements $txZ_0$ and $tx^{2^{m-1}+1}Z_0$ are $G$-conjugate but $x^{2^{m-1}}Z_0$ is not $G$-conjugate to these two. For both $i=1,2$, we have $\Ker(\Phi) = \Inn(R_{i}) \cong C_2 \times C_2$. Thus 
\[
	\Aut_{\overline{\mathcal{F}}}(R_1) \cong \Sigma_3,\qquad \Aut_{\overline{\mathcal{F}}}(R_2) \cong C_2.
\]
Therefore we get the same quiver as in Proposition~\ref{T:FusionExample}.

\subsection{Proof of Proposition \ref{T:SchurExample}} \label{S:SchurProof}

Let $B$ be the set of all upper triangular matrices in $G$. For $u \in \mathbb{F}_q$, set
\[
	[ u ] :=
	\begin{pmatrix}
		1 & u \\
		0 & 1
	\end{pmatrix}.
\]
Also set
\[
	t := 
	\begin{pmatrix}
		\epsilon & 0 \\
		0 & 1
	\end{pmatrix},
	\qquad\qquad
	w := 
	\begin{pmatrix}
		0 & 1 \\
		1 & 0
	\end{pmatrix}
\]
where $\epsilon$ is a generator of the multiplicative group $\mathbb{F}_q^{\times}$. Then we have
\[
	G / B = \{\, B, wB, [\epsilon^{i}]wB \,\}_{ 1 \leq i \leq q - 1}.
\]

Let 
\[
	kG \rightarrow \End_k (k[G/B])
\]
be the $k$-algebra homomorphism of Theorem \ref{T:SchurAlternative} and denote its image by $S$.  This map is the $k$-linear extension of the group homomorphism 
\[
	\psi : G \rightarrow \Sigma_{G/B} \hookrightarrow \GL_{k}(k [G/B])
\]
where the first homomorphism sends $g \in G$ to the permutation of $G/B$ induced by left multiplication by $g$ and the second inclusion sends permutations of $G/B$ to corresponding permutation matrices. Observe that the following correspondence
\begin{equation*}
\begin{array}{cccccc}
	B & wB & [ \epsilon ]wB & [ \epsilon^{2} ]wB & \cdots & [ \epsilon^{q-1} ]wB \\
	\updownarrow & \updownarrow & \updownarrow & \updownarrow & \cdots & \updownarrow \\
	\begin{bmatrix}
		1 \\
		0
	\end{bmatrix} &
	\begin{bmatrix}
		0 \\
		1
	\end{bmatrix} &
	\begin{bmatrix}
		\epsilon \\
		1
	\end{bmatrix} &
	\begin{bmatrix}
		\epsilon^{2} \\
		1
	\end{bmatrix} & \cdots &
	\begin{bmatrix}
		\epsilon^{q-1} \\
		1
	\end{bmatrix}
\end{array}
\end{equation*}
respects the $G$-action on $G/B$ by left multiplication and the natural $G$-action on the projective line over $\mathbb{F}_q$, where $\begin{bmatrix} u \\ v \end{bmatrix}$ denotes the image of $\begin{pmatrix} u \\ v \end{pmatrix}$ in the projective line. Denote above elements by $v_1$, $v_2$, $\ldots$ ,$v_{q+1}$, respectively, and write $V = k[G/B] = k v_1 \oplus k v_2 \oplus \dots \oplus k v_{q+1}$. Then $\psi$ factors through
\[
	\PGL_2(q) \cong G / Z(G) \hookrightarrow \GL_{k}(V),
\]
and hence
\[
	S = \Img ( k\PGL_2(q) \rightarrow \End_k (V) ).
\]

$V$ is a $(q+1)$-dimensional $S$-module with the natural $S$-action. Now we find its composition series. First of all, $V$ has an obvious $1$-dimensional simple $S$-submodule
\[
	V_1 = k ( v_1 + v_2 + \cdots + v_{q+1} ).
\]
Let us denote the elements of the quotient module $V / V_1$ as
\[
	[ \lambda_1, \lambda_2, \ldots, \lambda_{q+1} ] 
	:= \lambda_1 v_1 + \lambda_2 v_2 + \ldots + \lambda_{q+1} v_{q+1} + V_1
\]
with $\lambda_{i} \in k$. Then the $(q-1)$-dimensional $S$-submodule $V_2$ of $V / V_1$ given by
\[
	V_2 = \{\, [ \lambda_1, \lambda_2, \ldots, \lambda_{q+1} ] 
	\mid \lambda_1 + \lambda_2 + \ldots + \lambda_{q+1} = 0 \,\}
\]
is also simple because $\PGL_2(q)$ acts $3$-transitively on $\{\, v_1, v_2, \dots, v_{q+1} \,\}$.(See~\cite[Table 1]{Mortimer1980})  Let $W$ be the inverse image in $V$ of $V_2$.  Observe that $V$, $W$ are uniserial $S$-modules with composition series $(V_1, V_2, V_1)$, $(V_2, V_1)$, respectively.  In particular, both $V$ and $W$ are indecomposable.

It is well known that $V = k[G/B]$ is a projective $S$-module and that there are exactly two simple $S$-modules up to isomorphism.  Then, since $S$ is quasi-hereditary, it follows from the composition seris of $V$ that the standard modules for $V_1$ and $V_2$ are $V_1$ and $W$, respectively, and $W$ is also projective.  Therefore we conclude that $S$, and hence the $q$-Schur algebra $\mathcal{S}_2(q)$, is Morita equivalent to the path algebra of the quiver given in Proposition~\ref{T:SchurExample}.

\section{A remark on a canonical bijection of simple modules} \label{S:RemarkNaturalBijection}

Let $k$ be an algebraically closed field of characteristic $2$ and let $q$ be an odd prime power.  Let $b_0$ be the principal $2$-block of $G = \GL_n(q)$.  The algebra homomorphism in Theorem~\ref{T:SchurAlternative} restricts to the  surjective algebra homomorphism
\[
	kGb_0 \twoheadrightarrow S
\]
where $S$ is a $k$-algebra which is Morita equivalent to the $q$-Schur algebra $\mathcal{S}_n(q)$.  On the other hand, in Theorem~\ref{T:Main} we showed that there is another surjective algebra homomorphism
\[
	S_0 \twoheadrightarrow T_0
\]
where $S_0$ and $T_0$ are, respectively, the basic algebras of the $q$-Schur algebra $\mathcal{S}_n(q)$ and the weighted fusion category algebra $\overline{\F}(b_0)$ when $n=2$.  Combining these two surjective algebra homomorphisms, we see that simple $\overline{\F}(b_0)$-modules can be viewed as simple $kGb_0$-modules when $n=2$.  Since we have
\begin{equation*} \label{E:AWC}
	l(\overline{\F}(b_0)) = \text{number of partitions of } n = l(kGb_0) 
\end{equation*}
for $n = 2$ (in fact, for every $n$ by An~\cite{An1992}), we get a canonical bijection between simple $kGb_0$-modules and simple $\overline{\F}(b_0)$-modules in this case.  But, as mentioned in the remark after Theorem~\ref{T:FusionProperty}, there is a canonical bijection between the set of isomorphism classes of simple $\overline{\F}(b_0)$-modules and the set of conjugacy classes of $b_0$-weights.  Thus we get a canonical bijection between the set of isomorphism classes of simple $kGb_0$-modules and the set of conjugacy classes of $b_0$-weights when $n=2$.

\begin{ack}
We would like to thank Markus Linckelmann for his valuable suggestions and comments throughout writing this paper.  Also we thank Radha Kessar for her comment on the canonical bijection and Will Turner for pointing us to the work of Erdmann and Nakano~\cite{ErdmannNakano2001}.
\end{ack}

\bibliographystyle{amsplain}
\bibliography{mybib}
\end{document}